\input amstex
\input amsppt.sty
\magnification\magstep1
\input epsf
\def\ni\noindent
\def\sbs{\subset}

\def\as{\operatorname{asdim}}
\def\asdim{\operatorname{asdim}}
\def\diam{\operatorname{diam}}

\def\R{\text{\bf R}}

\def\Q{\text{\bf Q}}
\def\Z{\text{\bf Z}}

\def\N{\text{\bf N}}

\def\g{\gamma}
\def\G{\Gamma}

\def\mnorm#1{\| #1 \|}
\def\grpgen#1{ \langle #1 \rangle}
\def\Ucal{\Cal U}
\def\ds{\displaystyle}
\def\one{1_{\Gamma}}

\hoffset= 0.0in
\voffset= 0.0in
\hsize=32pc
\vsize=38pc
\baselineskip=24pt
\NoBlackBoxes
\topmatter
\author
J. Smith
\endauthor

\title
On asymptotic dimension of countable abelian groups
\endtitle
\abstract
We compute the asymptotic dimension of the rationals given with an
invariant proper metric. Also we show that a countable
torsion abelian group taken with an invariant proper metric has
asymptotic dimension zero.
\endabstract

\thanks
\endthanks

\address University of Florida, Department of Mathematics, P.O.~Box~118105,
358 Little Hall, Gainesville, FL 32611-8105, USA
\endaddress


\email  justins\@math.ufl.edu
\endemail

\endtopmatter

\document
\head \S1 Introduction \endhead

Gromov introduced the notion of asymptotic dimension as an invariant of finitely
generated discrete groups \cite{Gr}. This invariant was studied in a sequel of
papers \cite{B-D1},\cite{B-D2},\cite{B-D3},\cite{B-D-K},\cite{Ji} and others.
The notion of asymptotic dimension can be extended to the class of all
countable groups and most of the results for finitely generated groups
are valid for countable groups \cite{D-S}. To define asymptotic dimension for
a general countable group one should consider a left invariant proper
metric on it. It turns out that such metrics alway exists and any two
such metrics on a group $\Gamma$ are coarsely equivalent, i. e., they lead
to the same number $\asdim\Gamma$.

Even for very familiar infinitely generated countable groups, like
the group of rationals $\Q$, an invariant proper metric turns them into a
quite complicated geometrical object. To give an idea, we notice that
an invariant metric on $\Z[\frac{1}{2}]\subset\Q$ can be defined as
induced from the metric graph obtained by gluing together infinitely many
isosceles triangles with sides $2^n-2^n-2^{n-1}$, $n\in \N$, by the following
rule. First we glue triangles with sides 2-2-1 to all intervals $[n,n+1]\subset\R$,
$n\in\Z$ and mark their free vertices by the averages of the endpoints,
$\frac{2n+1}{2}$ (see Figure 1).

\medskip
\epsfysize=1.5in
\centerline{\epsffile{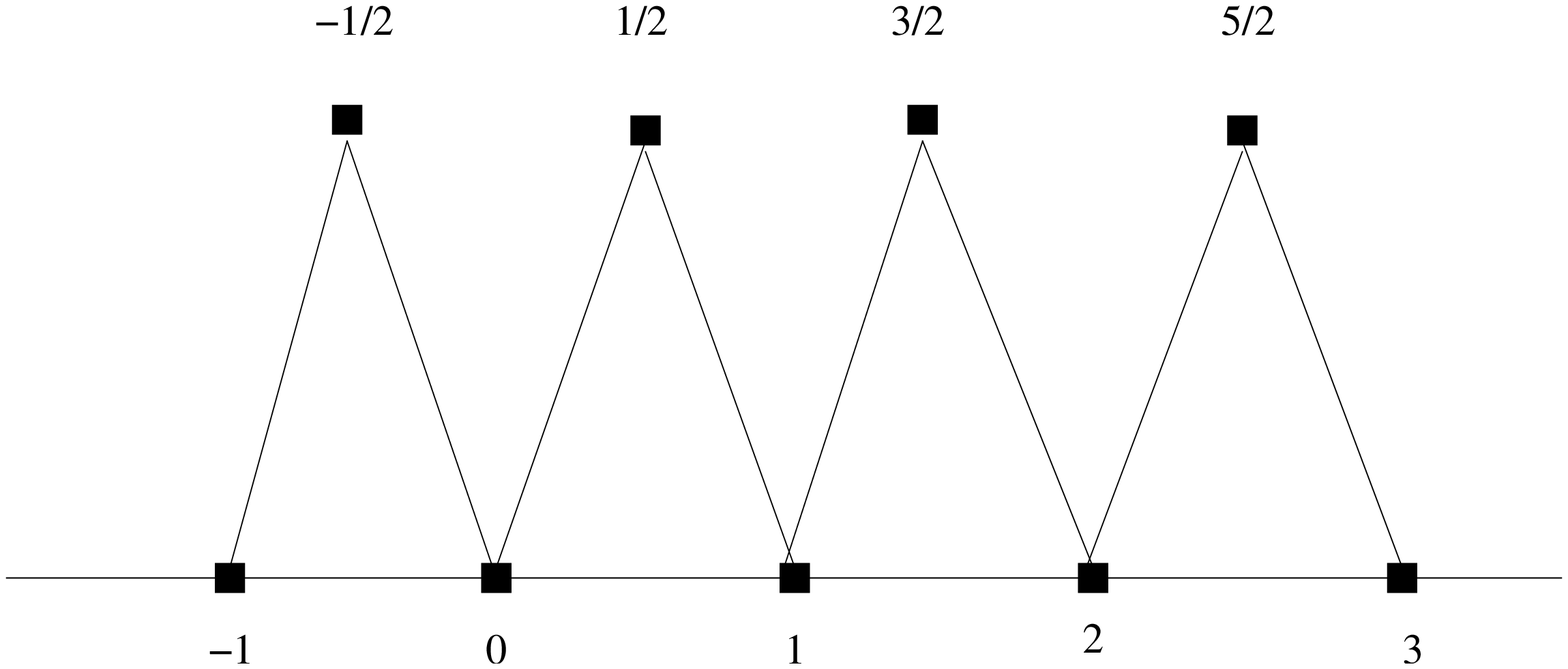}}
\smallskip
\centerline{Figure 1}
\medskip

Then to every edge of length 2 we glue a triangle with sides 4-4-2
and mark its free vertex by the average of the base and so on.
Then $\Z[\frac{1}{2}]$ is identified with the set of vertices of
this graph. 

The corresponding picture for $\Q$ is more
complicated. Nevertheless, in \S 3 we compute that $\asdim\Q=1$.
In particular, $\asdim\Z[\frac{1}{2}]=1$.

In the case of finitely generated groups, if $\asdim\Gamma=0$, then the group
$\Gamma$ is finite. This is not true for countable groups. In \S2 of the paper
we give a criterion for a group to have asymptotic dimension 0. As a corollary
we obtain that all torsion countable abelian groups have asymptotic dimension
0.

In \S 4 we show that the asymptotic dimension of the rationals taken with
$p$-adic norm is zero. Since $p$-adic norm is not proper on $\Q$ this cannot be
done by the criterion of \S 2.

\

{\bf Preliminaries.} 

The asymptotic dimension is defined for metric spaces.

DEFINITION~\cite{Gr}.  {\it We say that a metric space $X$ has asymptotic
dimension $\leq n$ if,
for every $d>0$, there is an $R$ and $n+1$ $d-$disjoint, $R-$bounded families
$\Ucal_0, \Ucal_1, \ldots, \Ucal_n$ of subsets of $X$ such that
$\cup \Ucal_i$ is a cover of $X$. }

We say that a family $\Ucal$ of subsets of $X$ is {\it $R-$bounded} if
$\sup \{ \diam{U} | U \in \Ucal \} \leq R$.  Also, $\Ucal$ is said to be
{\it $d-$disjoint}
if $d(x,y) > d$ whenever $x\in U$, $y \in V$, $U \in \Ucal$,
$V \in \Ucal$, and $U \neq V$.

The notion $\as$ is a coarse invariant (see \cite{Roe}).

Let $f: X \rightarrow Y$ be a map between metric spaces.  If, for
each $R>0$, there is an
$S>0$ such that $d(f(x),f(y)) <S$ whenever $d(x,y) < R$,
then we say that $f$ is {\it bornologous}.  If the
preimage of each bounded subset of $Y$ is a bounded subset of $X$,
then we say that $f$ is
{\it metrically proper}.  A map is said to be {\it coarse} if it is both
metrically proper
and bornologous.  Also, given two maps
$f,f' : X \rightarrow Y$, where $X$ is a set and $Y$ is a metric
space, then we say that
$f$ and $f'$ are close if $\sup_{x \in X} d(f(x),f'(x)) < \infty$.
We say that a metric space
$(X,d)$ is {\it proper} if closed, bounded sets are compact.

DEFINITION~\cite{Roe}. {\it Suppose $f: X \rightarrow Y$ is a coarse map
between metric spaces.  $f$ is a coarse equivalence
if there is a coarse map $g:Y \rightarrow X$ such that $f \circ g$ is close to the identity function on
$Y$ and $g \circ f$ is close to the identity function on $X$.}

 Let $G$ be a group.  A map $\mnorm{\cdot} :G \rightarrow [0, \infty) $
is said to be a {\it norm} on $G$ if $\mnorm{x}=0$ if and only if $x=1_{G}$,
$\mnorm{x^{-1}}=\mnorm{x}$ for all $x\in G$, and
$\mnorm{xy}\leq \mnorm{x}+\mnorm{y}$ for all $x, y\in G$.

Given a norm $\mnorm{\cdot}$ on $G$, define
$d:G \times G \rightarrow [0, \infty)$ by $d(x,y)=\mnorm{x^{-1}y}$, where
$x, y\in G$.  It is easy to verify that $d$ is a left invariant metric.  We say that a
norm on $G$ is {\it proper} if it has the property that for each $R>0$, there
are only finitely many $x\in G$ such that $\mnorm{x}\leq R$.  Then $d$ will induce the
discrete topology on $\Gamma$ and $d$ will be a proper metric.

Let $G$ be a finitely generated group with finite generating set $S$.
We define
$\mnorm{x} = \inf \{ n | x= \g_1 \g_2 \cdots \g_n, \g_i \in S\cup S^{-1} \}$.
This can be shown to be a norm on
the group $G$.  The metric induced by this norm is known as the {\it word metric}
on $G$ associated
with $S$.

Since the word metrics associated with any two finite generating sets of a finitely generated group are
coarsely equivalent (even
quasi-isometric), the asymptotic dimension is a group invariant.
Below we show that every countable group admits a left invariant proper metric.
In view of Proposition 1, one can extend the invariant $\as$ to
all countable groups (not necessarily finitely generated).

For countable groups, note that a left invariant, proper metric induces the discrete topology.  To see this, 
observe that for a proper metric, the group with this metric is complete.  To get a contradiction, suppose 
that there are no isolated points; then, as a consequence of the Baire category theorem, the group is 
not countable, a contradiction.  Thus, the group contains an isolated point.  As the metric is left invariant, 
left multiplication by a fixed element is an isometry and hence a homeomorphism.  This implies that 
every point in the group is isolated, so that the metric induces the discrete topology.  Thus, we have 
that a left invariant metric on a countable group is proper if and only if each 
bounded set is finite.  
\proclaim{Proposition 1}
For a countable group, any two left invariant, proper metrics are coarsely
equivalent.
\endproclaim
\demo{Proof}
Let $G$ be a group with left invariant, proper metrics $d$ and $d'$.  First, we show 
$id: (X,d) \rightarrow (X,d')$ is bornologous.  Let $R>0$ be given.  
Let $B(R,d)=\{ g\in G | d(g,1) \leq R \}$.  Since $d$ is proper,
$B(R,d)$ is finite.  Thus, there is an $S>0$ such that $B(R,d) \subset B(S,d')$.  So if $d(x,y)\leq R$, then
$d(1 ,x^{-1}y) \leq R$ since $d$ is left invariant.  Thus, $x^{-1}y \in B(R,d)$, and 
so $x^{-1}y \in B(S,d')$, or $d'(1 ,x^{-1}y) \leq S$.  Hence $d'(x,y) \leq S$.  So 
$id$ is bornologous.  By a similar argument, $id^{-1}$ is
bornologous.  This shows that $id$ and $id^{-1}$ are proper.  So $id$ is a coarse equivalence.
\qed
\enddemo

DEFINITION.  Let $\Gamma$ be a countable group.  Let $S$ be a generating set (possibly infinite)
for $\Gamma$.  A weight function $w:S\rightarrow [0, \infty)$ on $S$ is a function
such that the following properties hold:
\roster
\item{} if $w(s)=0$, then $1_{\Gamma} \in S$ and $s=1_{\Gamma}$,
\item{} $w(s)=w(s^{-1})$ whenever $s,s^{-1} \in S$, and
\item{} for each $N\in \N$, $w^{-1}[0,N]$ is a finite set.
\endroster

The third property says that $w$ is a proper map, where $\G$ has the discrete topology.
Also, this property can essentially be viewed as the requirement that
$\ds \lim w = \infty$.

It is not hard to see that for any countable group $\Gamma$, there is a weight function.
In fact, for any generating set $S$, there is a weight function with domain $S$.  Also, a weight
function $w:S\rightarrow [0, \infty )$ can be extended to a weight
function on $S \cup S^{-1}$ (or $S \cup S^{-1} \cup {1_{\Gamma}}$).

\proclaim{Theorem 1} A weight function on the countable group $\Gamma$ induces 
a proper norm $\mnorm{\cdot}$, and so
a weight function induces a left invariant, proper metric $d$. \endproclaim
\demo{Proof}  Given a weight function $w:S \rightarrow [0, \infty)$, where $S$ is a generating set
for the countable group $\Gamma$, define
$\ds \mnorm{x}=\inf{ \{\sum_{i=1}^{n}w(s_i)| x=s_1^{\epsilon_1}s_2^{\epsilon_2}\cdots %
s_n^{\epsilon_n}, s_i\in S, \epsilon_i\in \{ \pm 1 \} \} }$.  Note that if we view $\one$ as an empty product,
$\mnorm{\one} = 0$.  The proof that $\mnorm{\cdot}$ is a norm is left to the reader.

Let $R>0$ be given.  Let $r$ be a nonzero value that the weight function assumes (otherwise,
the weight function is always zero, in which case $\Gamma$ is trivial, and so the theorem
holds).  So $\{ s\in S | 0 < w(s)\leq r \}$ is nonempty and finite by definition.  Thus,
there is a $t\in S$ such that $w(t)= \min\{ w(s) | s \in S, 0<w(s)\leq r \}$.  It
is immediate that $0<w(t)\leq w(s)$ for all $s\in S \setminus \one$.  Now,
suppose $x$ is such that $\mnorm{x} \leq R$ and $x \neq \one$.  Then $\mnorm{x} < R+1$.  So there are
$s_1, s_2, \ldots, s_n \in S$, and $\epsilon_1, \epsilon_2, \ldots, \epsilon_n \in \{ \pm 1 \}$
such that $x=s_1^{\epsilon_1}s_2^{\epsilon_2}\cdots s_n^{\epsilon_n}$ and
$\sum_{i=1}^{n}w(s_i) < R+1$.  Further, we may assume that $s_i \neq \one$ for each $i$.
Thus, $s_i \in \{ s\in S | w(s)\leq R+1 \}$ for all $i$.  Also,
$R+1> \sum_{i=1}^{n}w(s_i) \geq n w(t)$, so that $n< (R+1)/w(t)$.  Thus,
$x$ is an element of $\{ t_1^{\delta_1}t_2^{\delta_2}\cdots t_m^{\delta_m} | t_i \in S \setminus \one,
w(t_i) \leq R+1, \delta_i \in \{ \pm 1 \}, m< (R+1)/w(t) \}$, a finite set.  This shows that
$\{ x \mid \mnorm{x}\leq R \}$ is finite.
\qed
\enddemo

Note that the infimum in the definition of $\mnorm{x}$ is actually a minimum.
To see this, simply modify the argument in the last paragraph to show that the set of elements of
$\{\sum_{i=1}^{n}w(s_i)| x=s_1^{\epsilon_1}s_2^{\epsilon_2}\cdots s_n^{\epsilon_n}, %
s_i\in S, \epsilon_i\in \{ \pm 1 \} \}$ less than $\mnorm{x}+1$ is a finite set.

\head \S2 Groups with asymptotic dimension 0 \endhead

The following theorem gives a necessary and sufficient for a
countable group to have asymptotic dimension zero.  This condition
relies only on the algebraic structure of the group.

\proclaim{Theorem 2} Let $G$ be a countable group.  Then $\asdim{G}=0$ if and only if every
finitely generated subgroup of $G$ is finite.
\endproclaim
\demo{Proof} 
Let $w:G \rightarrow [0,\infty)$ be a weight function on the generating set $G$.  Let
$\mnorm{\cdot}$ and $d$ be the induced norm and metric, respectively.

First suppose that $\asdim{G}=0$.  Let $T\subset G$ be a finite set.  Take
$d> \max_{g\in T}\mnorm{g}$.  As $\asdim{G}=0$, there is a $d-$disjoint, uniformly bounded
cover $\Cal{U}$ of $G$.  Choose $U \in \Ucal$  with $1 \in U$.  We will show that
$\grpgen{T} \subset U$.  To do this, we will show by induction that every product of $k$ ($k\geq 0$)
elements of $T \cup T^{-1}$ lies in $U$.  This is true for $k=0$, as $1\in U$.  Now suppose
it is true for $k-1$, $k\geq 1$.  Consider $x= t_1^{\epsilon_1}t_2^{\epsilon_2}\cdots t_k^{\epsilon_k}$,
where $t_i \in T$ and $\epsilon_i \in \{ \pm 1 \}$.  Set
$y= t_1^{\epsilon_1}t_2^{\epsilon_2}\cdots t_{k-1}^{\epsilon_{k-1}}$.  By the induction assumption,
$y \in U$.  Since $d(y,x) = \mnorm{y^{-1}x}= \mnorm{t_k^{\epsilon_k}}=\mnorm{t_k}< d$, and because
$\Ucal$ is a $d-$disjoint cover, we must have $x \in U$.  Thus, each product of $k$ elements of
$T \cup T^{-1}$ lies in $U$.  Therefore, $\grpgen{T} \subset U$.  As $\Ucal$ is uniformly bounded,
$U$ is bounded, and so $U$ and $\grpgen{T}$ are finite.

Conversely, suppose every finitely generated subgroup of $G$ is finite.  Let $d>0$ be given.  Define
$T= \{ s\in G | w(s) < d \}$ and $H = \grpgen{T}$.  By definition of weight function, $T$ is finite.
By our assumption, $H$ is finite as well.  Let $\Ucal = \{ gH | g \in G \}$ be the collection of
left cosets.  So $\Ucal$ is a uniformly bounded cover, as left multiplication by a fixed element is
an isometry of $G$.  Further, suppose $gH \neq hH$.  Let $x \in gH$ and
$y \in hH$.  So $xH=gH$ and $yH=hH$.  As $gH \neq hH$, $xH \neq yH$, and so $y^{-1}x \notin H$.  Hence
$y^{-1}x$ cannot be written as a product of elements of $T \cup T^{-1}$.  So if we take
$s_i \in G $ such that 
$y^{-1}x=s_1 s_2 \cdots s_n$ and $\mnorm{y^{-1}x}= \sum w(s_i)$,
then there is a $j$ such that $s_j \notin T$.  Hence $w(s_j) \geq d$, and so
$d(y,x) = \mnorm{y^{-1}x} \geq d$.  Therefore $\Ucal$ is a $d-$disjoint, uniformly bounded cover.
Since $d>0$ was arbitrary, $\asdim{G}=0$.  This completes the proof.
\qed
\enddemo
The following corollaries are immediate consequences.

\proclaim{Corollary 1} Let $G$ be a finitely generated group.  Then $\asdim{G}=0$ iff $G$ is a finite
group.  \endproclaim

\proclaim{Corollary 2} Let $G$ be a countable abelian group.
Then $\asdim{G}=0$ if and only if $G$ is a torsion group.  \endproclaim

REMARK.  The last corollary shows that  $\oplus_i\Z_{m_i}$, $\Q / \Z$,
and $\Z_{p^{\infty}}=\lim_{\to}\Z_{p^k}$ all have asymptotic dimension 0.

The next theorem states that the epimorphic image of a zero dimensional countable group is zero
dimensional. This is not true for one dimensional groups. Moreover, every
countable group is an epimorphic image of a free group which is one dimensional.

\proclaim{Theorem 3} Let $\phi : G \rightarrow H$ be an epimorphism of countable groups.  If
$\asdim{G}=0$, then $\asdim{H}=0$.  \endproclaim

\demo{Proof}   We will show that every finitely generated subgroup of $H$ is finite.  Let
$T$ be a finite subset of $H$.  Since $\phi$ is an epimorphism, for each $t \in T$ we can find
a $g_t \in G$ for which $\phi (g_t) =t$.  Let $S = \{ g_t | t \in T \}$.  As $S$ is finite, and
since $\asdim{G}=0$, we have that $\grpgen{S}$ is finite by the theorem.  Thus,
$\grpgen{T} = \grpgen{\phi (S)} = \phi( \grpgen{S} )$ is a finite set.  By the theorem, $\asdim{H}=0$.
\qed
\enddemo

\head \S3 Asymptotic dimension of the rationals \endhead

We will now show that $\ds \asdim{\Q}=1$.  Once more we note that  we are
not computing
the dimension of $\Q$ with the Euclidean metric here, but rather with
a {\bf proper},
left-invariant metric.

\proclaim{Theorem 4} $\as\Q=1$. \endproclaim
\demo{Proof}
First, we will show that $\ds \asdim{\Q \cap [0, 1)}=0$.  We
will then use this to prove the result.

On $\Q$, define $\mnorm{m/n}_{\Q}=|m/n|+\ln(n)$ when $m$ and $n$ are relatively prime integers and $n$
is positive (when $m$ and $n$ have these properties, we will say that $m/n$ is in standard form).
Here $|\cdot |$ denotes the usual absolute value.  It is easy to
show that $\mnorm{\cdot}_{\Q}$ is a proper norm on $\Q$.  This norm induces a metric
$d_{\Q}$ in the usual way.  Also, it is not hard to show
that, on $\Q \cap (-1, 1)$, $\ln{n} \leq \mnorm{m/n}_{\Q} \leq 3\ln{n}$ whenever
$m/n$ is in standard form.

Let $\G=\Q / \Z$.  Let $p:\Q \rightarrow \G$ be the usual projection map.  $p$ is a surjective homomorphism and
$\ker p=\Z$.  As the topology induced by $\mnorm{\cdot}_{\Q}$ is discrete,
$\Z$ is closed in $\Q$.  Thus, this norm induces a norm on $\Q / \Z$, given by
$$\mnorm{\overline{x}}=\inf{ \{ \mnorm{x+m'}_{\Q}\ |\ m'\in \Z \} }.$$  This is a proper norm
and hence its associated metric $d$ will be a left-invariant, proper metric on $\Q / \Z$.

Define $i:\Q / \Z \rightarrow \Q\cap[0,1)$ as follows:  For $r\in \Q$, there is a
unique $r'\in \Q \cap [0,1)$
such that $\overline{r}=\overline{r'}$; set $i(\overline{r})=r'$.  It is easy to see that $i=p|_{\Q\cap [0,1)}^{-1}$.

We will show that $i$ is a coarse equivalence.  First, $\mnorm{p(r)}=\mnorm{\overline{r}}\leq \mnorm{r}_{\Q}$,
so $p$ and hence $p|_{\Q\cap [0,1)}$ is bornologous.  Further, suppose $\frac{m}{n}\in \Q\cap (-1,1)$,
$(m,n)=1$, $m, n \in \Z$, and $n>0$.  For $m'\in \Z$, we know that $(m+nm', n)=1$.
Thus,
$$
\mnorm{\frac{m}{n}+m'}_{\Q}=\mnorm{\frac{m+nm'}{n}}_{\Q}= 
|\frac{m}{n}+m'|+\ln(n)\geq \ln(n).
$$  
Also, since $-1<\frac{m}{n}<1$, we have that
$\ln(n)\geq \frac{1}{3}\mnorm{\frac{m}{n}}_{\Q}$.  Hence
$$
\mnorm{\overline{\frac{m}{n}}}=\inf{ \{ \mnorm{\frac{m}{n}+m'}_{\Q}\ |\ m'\in \Z \} } \geq \ln(n)\geq %
\frac{1}{3}\mnorm{\frac{m}{n}}_{\Q}.$$  Since each $r \in \Q\cap (-1,1)$ can be
expressed in standard form,
$\mnorm{r}_{\Q} \leq 3\mnorm{\overline{r}}$.  So
for $r, s \in\Q\cap [0,1)$, we have $s -r \in \Q\cap (-1,1)$ and so
$$d_{\Q}(i(\overline{r}), i(\overline{s}))=%
\mnorm{-i(\overline{r}) +i(\overline{s})}_{\Q}=%
\mnorm{s-r}_{\Q}\leq 3\mnorm{\overline{s-r}}%
= 3\mnorm{\overline{s}-\overline{r}}=3d(\overline{r},\overline{s}).$$  Since for each $r \in \Q$, there
is a $r' \in \Q\cap [0,1)$ such that $\overline{r}=\overline{r'}$, we have that
$d_{\Q}(i(\overline{r}), i(\overline{s})) \leq d(\overline{r},\overline{s})$ for $r,s \in \Q$.  Thus,
$i$ is bornologous.  As
$p|_{\Q\cap [0,1)}$ and $i$ are inverses, each is proper and $i$ is a coarse equivalence of
$\Q / \Z$ and $\Q \cap [0,1)$.  By Corollary 2,
$\asdim{ \Q \cap [0,1)}=\asdim{\Q / \Z}=0$.

We will now complete the proof that $\asdim{\Q}=1$.  Let $d>0$ be given.  Since
$\{ x\in\Q\ |\ \mnorm{x}_{\Q}\leq d\}$ is a finite set, there is an $R\in \Z_{+}$ such that
$\mnorm{x}_{\Q} \leq d$ implies $|x|< R$.  For $n\in \Z$, define $A_n=\Q \cap [nR, (n+1)R)$.
Notice $\asdim{A_0}=0$ by the finite union theorem of \cite{B-D1}.  So there is a $d-$disjoint,
$S-$bounded covering $\{ A_{0,k} | k=1,2,\ldots \}$ of $A_0$.  Since the map
$x \rightarrow x+nR$ ($n$ fixed) is an isometry $[0, R] \rightarrow [nR, (n+1)R]$, the
covering $\{A_{n,k} | k=1, 2, \ldots \}$ of $A_n$, where $A_{n,k}=nR + A_{0,k}$, is
$d-$disjoint and $S-$bounded.

Let $$\Cal{U}_0=\{ A_{n,k} | n \text{\ even}, k=1,2,\ldots \} \ \ \text{and}\ \ \
\Cal{U}_1=\{ A_{n,k} | n \text{\ odd}, k=1,2,\ldots \}.$$
Note that $\Cal{U}_0 \cup \Cal{U}_1$
is a cover of $\Q$, and $\Cal{U}_0$ and $\Cal{U}_1$ are both $S-$bounded.  We will
now show that $\Cal{U}_0$ is $d-$disjoint.  Consider $A_{n,k}$ and $A_{n',k'}$, where
$(n,k)\neq (n',k')$ and $n,n'$ are even.  Suppose first that $n\neq n'$.  Without loss of generality,
take $n<n'$.  Let $x\in A_{n,k}$ and $y\in A_{n',k'}$.  Then since
$A_{n,k}\subset A_n$ and $A_{n',k'} \subset A_{n'}$, we have
$nR\leq x< (n+1)R$ and $n'R\leq y< (n'+1)R$, and so
$y-x \geq (n'-n-1)R \geq R$.  Hence $|y-x| \geq R$.  But by our choice of $R$, this implies
$\mnorm{y-x}_{\Q} > d$.  Now we will consider the case when $n=n'$.  This forces $k\neq k'$.  By our
construction of the $A_{n,k}$, $A_{n,k}$ and $A_{n',k'}=A_{n,k'}$ are $d-$disjoint.

Similarly, $\Cal{U}_1$ is a $d-$disjoint family.  Since $d>0$ was arbitrary, $\asdim{\Q}\leq 1$.
Finally, since $d_{\Q}$ restricts to the Euclidean metric on $\Z$, we have $\asdim{\Q} \geq 1$.
Therefore, $\asdim{\Q}=1$.
\qed
\enddemo

\head \S4 Asymptotic dimension of the rationals with $p$-adic norm\endhead

We will now consider $\Q$ with the $p$-adic norm $\mnorm{\cdot}_p$.  Namely, if $m=p^a m'$ and
$n=p^b n'$, where $p$ divides neither $m'$ nor $n'$, then
$$\mnorm{\frac{m}{n}}_p=p^{b-a}=\frac{1}{p^{a-b}}.$$
Let $d_p$ denote the metric obtained from this norm.
Unlike the previous examples, $\Q$ with the metric $d_p$ is not proper.
To differentiate the dimension with respect to this metric from the one in Theorem 4, we
will always write $\asdim{(\Q, \mnorm{\cdot}_p)}$.
\proclaim{Theorem 5} $\as(\Q,\mnorm{\cdot}_p)=0$.
\endproclaim
\demo{Proof} Set $L=\{ \frac{m}{p^a} \mid a\in\N, m\in\Z \}$ and $X=L \cap [0,1)$.

We will now show that $N_1(X)=\Q$, where $N_1(X)= \{ y \in \Q \mid d_p(y,X) \leq 1 \}$.  
Suppose $r \in \Q$.  If $r=0$, then 
$r \in X \subset N_1 (X)$.  Now suppose $r \neq 0$.  Then there are $m',n' \in \Z \setminus 0$ 
such that $r = \frac{m'}{n'}$.  So $m'=p^a m$ and
$n'=p^b n$ for some $a,b \in\N$ and $m,n \in \Z \setminus 0$ such that $p$ divides 
neither $m$ nor $n$.  Thus, 
$\frac{m'}{n'}=p^{a-b}\frac{m}{n}$.  First consider the case when $a\geq b$. Then
$\mnorm{\frac{m'}{n'}}_p=\frac{1}{p^{a-b}}\leq 1$.  So $r =\frac{m'}{n'} \in N_1(0)\subset N_1(X)$.
Now suppose $a<b$.  Set $c=b-a$, so $\frac{m'}{n'}=\frac{m}{p^c n}$.  Since $(n, p^c)=1$,
$\overline{n} \in \Z_{p^c}$ is a generator.  Hence there exists an $\ell$ such that
$0\leq \ell < p^c$ and $\overline{m}=\ell \overline{n}$.  Therefore, $p^c\mid m-\ell n$.
Now take $d\geq 0$ such that $p^d \mid m-\ell n$ yet $p^{d+1}$ does not divide $m-\ell n$.  
So $d\geq c$.
Thus,
$$
\mnorm{\frac{m'}{n'}-\frac{\ell}{p^c}}_p=\mnorm{\frac{m}{p^c n}-\frac{\ell}{p^c}}_p=%
\mnorm{\frac{m-\ell n}{n p^c}}_p=\frac{1}{p^{d-c}}\leq 1.
$$
As $\frac{\ell}{p^c}\in X$, $r = \frac{m'}{n'} \in N_1 (X)$.  Therefore, $\Q=N_1 (X)$.

This means that the inclusion $X \hookrightarrow \Q$ is a coarse equivalence (see \cite{Roe}), where $X$ has the
restricted metric $d_p$.
Thus, $\asdim{(X,d_p) }=\asdim{(\Q,\mnorm{\cdot}_p)}$.  Let $\mnorm{\cdot}_{\Q}$ be the
norm from theorem 4, and we consider its restriction to $X' = L \cap (-1,1)$.  
Let $r \in X'\setminus 0$.  So $r = \frac{m}{p^a}$ for some 
$a \geq 0$ and $m \in \Z \setminus 0$ such that $(m,p^a)=1$.  
Since we also have $-1<\frac{m}{p^a} < 1$, it follows 
that $p$ does not divide $m$.  Then, since $\frac{m}{p^a}$ is in standard form and
$\frac{m}{p^a} \in (-1,1)$, we have
$$
\ln{(p^a)} \leq \mnorm{\frac{m}{p^a}}_{\Q} \leq 3\ln{(p^a)}.
$$
But $p^a = \mnorm{\frac{m}{p^a}}_p$ since $m \neq 0$ and $p$ does not divide $m$.  So
$$\ln(\mnorm{\frac{m}{p^a}}_p) \leq \mnorm{\frac{m}{p^a}}_{\Q} \leq
3 \ln(\mnorm{\frac{m}{p^a}}_p).$$
Thus, $\ln(\mnorm{r}_p) \leq \mnorm{r}_{\Q} \leq 3 \ln(\mnorm{r}_p)$.  For
$x,y \in X$ such that $x \neq y$, we have $y-x \in X' \setminus 0$, and so
$$\ln( d_p(x,y)) \leq d_{\Q} (x,y) \leq 3 \ln(d_p (x,y)).$$  From this it
is immediate that $id:(X, d_{\Q}) \rightarrow (X,d_p)$ and its inverse are
bornologous, and so they are coarse equivalences as well.  From the results of Theorem 4,
$\asdim{(X, d_p)} = \asdim{(X,  d_{\Q})} \leq \asdim{(\Q \cap [0,1), d_{\Q})} = 0$.
Thus, $\asdim{(\Q, \mnorm{\cdot}_p)}=0$.
\qed
\enddemo

\enddocument